\documentclass[12pt]{amsart}
\usepackage{amssymb}
\usepackage{amscd}
\usepackage{amsmath}

\newcommand{\gb}{{\mathfrak b}}
\newcommand{\gog}{{\mathfrak g}}
\newcommand{\gh}{{\mathfrak h}}
\newcommand{\gl}{{\mathfrak l}}
\newcommand{\gm}{{\mathfrak m}}
\newcommand{\gn}{{\mathfrak n}}
\newcommand{\gp}{{\mathfrak p}}
\newcommand{\gs}{{\mathfrak s}}

\newcommand{\bB}{{\bf B}}
\newcommand{\bG}{{\bf G}}
\newcommand{\bL}{{\bf L}}
\newcommand{\bM}{{\bf M}}
\newcommand{\bP}{{\bf P}}

\newcommand{\al}{\alpha}
\newcommand{\be}{\beta}

\newcommand{\cal}{\mathcal} %% Needed for LaTeX2e
\newcommand{\Oscr}{{\cal O}}
\newcommand{\Vscr}{{\cal V}}
\newcommand{\Iscr}{{\cal I}}
\newcommand{\Lie}{{\rm Lie\,}}
\newcommand{\QED}{\par \hspace{15cm}$\blacksquare$ \par}
\newcommand{\pr}{^{\prime}}
\newcommand{\st}{\subset}
\newcommand{\Pf}{\noindent{\bf Proof.}\par\noindent}

\newcommand{\sr}{\scriptscriptstyle}
\newcommand{\ov}{\overline}
\newcommand{\rar}{\rightarrow}
\newcommand{\Spa}{{\rm span\,}}
\newtheorem*{theorem}{Theorem}
\newtheorem*{lemma}{Lemma}
\newtheorem*{prop}{Proposition}
\newtheorem*{conj}{Conjecture}
\begin{document}
\title[orbital varieties]
{\bf On varieties in an orbital variety closure in semisimple Lie algebras}
\author{Anna Melnikov}
\address{Department of Mathematics,
University of Haifa,
Haifa 31905, Israel}
\email{melnikov@math.haifa.ac.il}
\begin{abstract}
In this note we discuss the closure of an 
orbital variety as a union of varieties. We
show that if semisimple Lie algebra $\gog$ contains factors 
not of type $A_n$
then there are orbital varieties whose closure contains
components which are not Lagrangian. We show that the argument
does not work if all the factors are of type $A_n$ and
provide the facts supporting the conjecture claiming that if
all the factors of $\gog$ are of type $A_n$ then the
closure of an orbital variety is a union of orbital
varieties.
\end{abstract}
\maketitle

\section {\bf Introduction}
\subsection{}\label{1.1}
Let $\bG$ be a connected simply-connected complex algebraic group. Set $\gog=\Lie(\bG).$
Consider the co-adjoint action of $\bG$ on $\gog^*.$ Identify $\gog^*$ with $\gog$
through the Killing form. A $\bG$ orbit $\Oscr$ in $\gog$ is called nilpotent
if it consists of ad-nilpotent elements. 

Fix a triangular decomposition $\gog=\gn^-\oplus\gh\oplus\gn.$  
\par
Let $\Oscr$ be some nilpotent orbit. 
Consider an intersection $\Oscr\cap\gn.$ As it was shown by N. Spaltenstein \cite{Sp} and
R. Steinberg \cite{St}, this is an equidimentional variety of the dimension $0.5\dim \Oscr.$
Moreover, it was shown by A. Joseph \cite{J} that this is a Lagrangian subvariety of $\Oscr.$

An irreducible component of 
$\Oscr\cap\gn$
is called an orbital variety associated to $\Oscr.$ 

According to orbit method
philosophy, one would like to attach an 
irreducible representation of the enveloping algebra $U(\mathfrak g)$ to an orbital variety. This
should be a simple highest weight module.
The results of A. Joseph and T. A. Springer
provide  
a one to one correspondence
between the set of primitive ideals of $U(\mathfrak g)$ 
containing the augmentation ideal of its centre (thus, corresponding to
integral weights) 
and the set of orbital varieties in $\mathfrak g$ corresponding
to Lusztig's special orbits.
Thus, orbital varieties play a key role 
in the study of
primitive ideals in $U(\mathfrak g).$ The details can be found in
\cite{B-B}, \cite{J} and \cite{M1}.

\subsection{}\label{1.2}
The closure of a nilpotent orbit is a union
of nilpotent orbits. The combinatorial description of this union was given
by M. Gerstenhaber for $\gog=\gs\gl_n.$ Further H. Karft and C. Procesi
described this union for other simple Lie algebras.

Generalizing the results of N. Spaltenstein, D. Mertens \cite{M}
showed that $\ov{\Oscr\cap\gn}=\ov\Oscr\cap\gn.$ Thus,  $\ov{\Oscr\cap\gn}$ is a union
of intersections of $\gn$ with corresponding orbits defined by the results of Gerstenhaber 
and Kraft-Procesi.

The question is to give a description of an orbital variety closure in the spirit of
Gerstenhaber theory.
This question is much more involved than the question on a nilpotent orbit closure.
It has two components.
The first one , a purely geometrical component, is to describe the type of varieties
which  constitute the closure of an 
orbital variety. This question can be formulated as following. Let $\Vscr$ be an orbital
variety. Then its $\bG$-saturation $\Oscr_{\Vscr}$  is  a nilpotent orbit, to which $\Vscr$ is associated.
Let us take some nilpotent orbit  $\Oscr\ :\ \Oscr\st\ov\Oscr_{\Vscr}$ and consider $\ov\Vscr\cap \Oscr.$
As we show in \ref{3.1}, this intersection is always not empty. So, a natural task
is to describe the irreducible components of this intersection. Is this intersection equidimentional? 
Is this intersection Lagrangian?

Here we show  that if $\gog$ contains factors not of type $A_n,$ there exist
orbital varieties in $\gog$ such that the intersection mentioned above
is not Lagrangian. 

We explain why the same argument does not work if all factors
are of type $A_n.$
Moreover, the study of special cases in \cite{M2} and ~\cite{M3} shows that at least
for some special types of orbital varieties in $\gs\gl_n$ (that is $A_{n-1}$) the intersection 
of an orbital variety closure $\ov\Vscr$ with any nilpotent orbit in the closure of $\Oscr_\Vscr$ is 
Lagrangian. 
Together with computations in low rank cases (for $n\leq 6$) this supports 
\begin{conj} In $\gs\gl_n$ the closure of an orbital variety is a union of orbital
varieties.
\end{conj}

Note that in any case the straight generalization of Gerstenhaber theory cannot work
for orbital varieties. As it is shown in \cite[4.1]{M2}, even if $\Vscr$ is of the most 
simple type (that is, when $\ov\Vscr$ is a nilradical) in $\gs\gl_n$ and $\Oscr\st \ov\Oscr_{\Vscr}$ 
one has that in general $\ov{\ov\Vscr\cap\Oscr}\ne\ov\Vscr\cap\ov\Oscr.$

\subsection{}\label{1.3}
The second component of the description of an orbital variety closure is
to give a combinatorial algorithm describing all the orbital varieties included in
the closure of a given one. This is a very complex combinatorial question.
The only general description of an orbital variety is provided by Steinberg's
construction (cf. \ref{2.1}). It is given via surjection from the Weyl group
onto the set of orbital varieties. But even the description of the fibers 
of this map is highly non trivial outside of type $A_n.$ For the type $A_n$
the picture is much nicer and simpler. Here the fibers are described by
Robinson-Schensted procedure and the question can be formulated in
terms of partial ordering of Young tableaux. Partial results in this case
are provided in \cite{M1}, \cite{M2}, \cite{M3}.

\subsection{}\label{1.4}
The body of the paper consists of two sections. In section 2 we give all the essential
background to make this note self-contained. In section 3 we provide the results 
on the orbital variety closures.

\section{\bf Preliminaries}
\subsection{}\label{2.1}
The definition of an orbital variety does not provide any way to
construct it. The only general construction of an orbital variety
belongs to R. Steinberg \cite{St}. We explain it here in short. 

Let $\gog$ be any semisimple Lie algebra. Fix its triangular decomposition
$\gog=\gn^-\oplus\gh\oplus\gn.$ For any $X\in \gn$ put $\Oscr_X:=\{gXg^{-1}\ :\ g\in\bG\}$
to be its nilpotent orbit.

Let $\bB$ be the Borel subgroup of $\bG$ with $\Lie(\bB)=\gh\oplus \gn$ and let
$\bB$ act adjointly on $\gn.$ 

Let $R\st\gh^*$ denote the set of non-zero roots,  $R^+$ the set
of positive roots corresponding to $\gn$ in the triangular 
decomposition of $\gog,$ and $\Pi\st R^+$ the resulting set of 
simple roots. Let $X_\alpha$ denote the root subspace for $\alpha\in R.$ 
One has $\gn=\bigoplus\limits_{\alpha\in R^+}X_\alpha.$

Let $W$ be the Weyl group of $<\gn,\gh>.$
The action of $w\in W$ on root subspace $X_\alpha$ is defined (in a standard way)
by $w(X_\alpha)=X_{w(\alpha)}.$ Consider the following subspace of $\gn:$
$$\gn\cap^w\gn=\bigoplus\limits_{\{\alpha\in R^+\, |\, w^{-1}(\alpha)\in R^+\}}X_\alpha.$$
Consider $\ov{\bG(\gn\cap^w\gn)}.$ Since the number of orbits is finite, this is
a closure of the unique orbit which we denote by $\Oscr_w.$
By R. Steinberg \cite{St}, one has
\begin{theorem} For each $w\in W$ there exists an orbital variety 
$\Vscr$ and  for each orbital variety  
$\Vscr$ there exists $w\in W$ such that 
$$\Vscr=\ov{\bB (\gn\cap{^w}\gn)}\cap\Oscr_w.$$
\end{theorem}
In what follows we will denote $\Vscr_w:=\Vscr$ in that case.
Obviously, $\Vscr_w$ is associated to $\Oscr_w.$

\subsection{}\label{2.2}
For any $\alpha\in \Pi$ let $\bP_\alpha$ be the standard parabolic subgroup of $\bG$
such that $\Lie(\bP)=\gb\oplus X_{-\alpha}.$

Given $\Iscr\st \Pi,$ let ${\bP}_{\Iscr}$ denote the unique 
standard parabolic subgroup of $\bG$ such that 
$\bP_{-\alpha}\st \bP_\Iscr$ iff $\alpha\in\Iscr.$ Let ${\bM}_{\Iscr}$ be the unipotent 
radical of ${\bP}_{\Iscr}$ and ${\bL}_{\Iscr}$ a Levi factor.
Let $\gp_{\sr \Iscr},\ \gm_{\sr \Iscr},\ \gl_{\sr \Iscr}$ denote 
the corresponding Lie algebras.
Set ${\bB}_{\Iscr}:={\bB}\cap{\bL}_{\Iscr}$ and
$\gn_{\sr \Iscr}:=\gn\cap\gl_{\sr \Iscr}.$
 We have decompositions ${\bB}={\bM}_{\Iscr}\ltimes{\bB}_\Iscr$ and 
$\gn=\gn_{\sr \Iscr}\oplus\gm_{\sr \Iscr}.$ They define projections 
${\bB}\rar {\bB}_\Iscr$ and $\gn\rar\gn_{\sr \Iscr}$ which we denote by 
$\pi_{\sr \Iscr}.$

Set $W_\Iscr:=<s_{\al}\ :\ \al\in \Iscr>$ to be a 
parabolic subgroup of $W.$ Set
$D_\Iscr:=\{w\in W\ :\ w(\al)\in R^+\ \forall \ \al\in \Iscr\}.$ 
Set $R_\Iscr^+=R^+\cap \Spa(\Iscr).$
A well known result provides that 
each $w\in W$ has a unique expression of the 
form $w=w_{\sr \Iscr}d_{\sr \Iscr}$
where $d_{\sr \Iscr}\in D_\Iscr,\ w_{\sr \Iscr}\in W_\Iscr$ and $\ell(w)=\ell(w_{\sr \Iscr})+
\ell(d_{\sr \Iscr}).$  Moreover, for any $w\in W$ one has 
$$R_\Iscr^+\cap^w R^+=R^+_\Iscr\cap^{w_{\sr \Iscr}} R^+_\Iscr.$$
Thus, a decomposition $W=W_\Iscr\times D_\Iscr$   
defines a projection $\pi_{\sr \Iscr}:W\rar W_\Iscr.$ For $w\in W$ set $w_{\sr \Iscr}:=\pi_{\sr \Iscr}(w).$
This element can be regarded as an element of $W_\Iscr$ and as an element of $W.$

Let ${\Vscr}_{w_{\sr \Iscr}}$ be the corresponding orbital variety
in $\gog$ and ${\Vscr}^{\Iscr}_{w_{\sr \Iscr}}$ be the corresponding orbital variety
in $\gl_{\sr \Iscr}.$ As it is shown in \cite[4.1.1]{M1} all the projections  
are in correspondence on orbital varieties, namely
\begin{theorem} Let $\gog$ be a reductive Lie algebra. Let $\Iscr\st \Pi.$ 
For every orbital variety ${\Vscr}_w\st \gog$ one has $\pi_{\sr \Iscr}(\ov{\Vscr}_w)=
\ov{{\Vscr}^{\Iscr}_{w_\Iscr}}.$ 
\end{theorem}

\subsection{}\label{2.3}
In what follows we need also 
the notion of $\tau$-invariant.
Let $w$ be any element of $W.$ Set  $S(w):=R^+ \cap^w R^-=\{\al\in R^+\ :\ w^{\sr -1}(\al)\in R^-\}.$
Set $\tau(w)=\Pi\cap S(w).$

As it
can be seen immediately from Steinberg's construction for orbital variety closures, one has
\begin{prop} Let $w,y\in W.$
If $\ov\Vscr_w\st \ov\Vscr_y$ then $\tau(w)\supset \tau(y).$
\end{prop}
Note that as a trivial corollary we get that if $\Vscr_w=\Vscr_y$ then $\tau(w)=\tau(y).$
In other words, $\tau$ invariant is constant on an orbital variety. 

\section{\bf An orbital variety closure}

\subsection{}\label{3.1}
We begin with a simple
\begin{lemma} 
Let $\Oscr,\Oscr\pr\st\gog$ be two nilpotent orbits
such that  $\ov\Oscr\pr\subsetneq\ov\Oscr.$ Let $\Vscr$ be an 
orbital variety associated to $\Oscr.$ Then
$\ov\Vscr\cap\Oscr\pr\ne 
\emptyset$ and in particular
there exist an orbital variety $\Vscr\pr$ 
associated to $\Oscr\pr$ such that
$\ov\Vscr\cap\Vscr\pr\ne \emptyset.$
\end{lemma}
\Pf
Indeed since there exist $w\in W$ such that 
$\ov\Vscr=\ov{\bB(\gn\cap^w\gn)}$
 and 
since $\bG/\bB$ is projective we get
$$\Oscr\pr\st\ov\Oscr=\ov{\bG(\gn\cap^w \gn)}=
\bG(\ov{\bB(\gn\cap^w\gn})=\bG(\ov\Vscr)$$
which proves the first part. 
\par
Since $\Oscr\pr\cap\ov\Vscr=\Oscr\pr\cap\gn
\cap\ov\Vscr$ we get the existence of $\Vscr\pr.$
\QED
\subsection{}\label{3.2}
Our strategy is to show that if in $\gog$ not 
all the factors are of type
$A_n$, there exist nilpotent orbits $\Oscr_1,\Oscr_2$ such that 
$\ov\Oscr_2\subsetneq\ov\Oscr_1$ and there exists  $\Vscr_w$ associated 
to $\Oscr_1$ such that
for every $\Vscr_z$ associated to $\Oscr_2$ 
one has  $\tau(w)\not\st
\tau(z).$ 
Then on one hand by lemma \ref{3.1} there exist at least one
$\Vscr_z$ 
associated to $\Oscr_2$ such that
$\ov\Vscr_w\cap\Vscr_z\ne \emptyset.$
On the other hand by proposition \ref{2.3}  if
$\tau(z)\not\supset
\tau(w)$ then $\Vscr_z\not\st\ov\Vscr_w$
for every $\Vscr_z$ associated to $\Oscr_2.$
We get that $\ov\Vscr_w\cap \Oscr_2$ is a non empty variety
of dimension less than $0.5\dim\Oscr_2.$
\subsection{}\label{3.3}
Consider the algebras of type $B_2$ and $G_2.$ 
They are fully described in \cite{T} and we just follow
these computations.

{\bf $B_2:$}\ \ \ Let $s$ be a reflection 
corresponding to the 
long root $\be$ and $t$ be the reflection corresponding 
to the short root $\al .$ Consider $\Oscr_s$ and 
$\Oscr_{tst}.$ One has
$\ov\Oscr_s\supsetneq 
\Oscr_{tst}.$ 
Moreover, $\Oscr_{tst}\cap \gn$ is irreducible so 
$\Vscr_{tst}$ is the unique orbital 
variety associated to $\Oscr_{tst}.$ 
Consider $\Vscr_s.$ 
If its closure is a union of orbital varieties 
then by lemma \ref{3.1}  $\Vscr_{tst}$ must be included
in it. But $\tau(s)=\{\be\}$ and 
$\tau(tst)=\{\al\}$ so this inclusion contradicts 
proposition \ref{2.3}.

In what follows we will need also the following fact 
about these orbits: there is no intermediate nilpotent orbit
between $\Oscr_s$ and $\Oscr_{tst}$, that is if $\Oscr\pr$ is 
such that $\ov\Oscr_{tst}\subset\ov{\Oscr\pr}\subset\ov\Oscr_{s}$
then $\Oscr\pr=\Oscr_{tst}$ or $\Oscr\pr=\Oscr_s.$

{\bf $G_2:$}\ \ \ In that case the picture is very 
similar to $B_2.$ Let $s$ be a reflection 
corresponding to the 
long root $\be$ and $t$ be the reflection corresponding 
to the short root $\al.$
Once more $\ov\Oscr_s\supsetneq \ov\Oscr_{tstst}$ and
$\Oscr_{tstst}\cap \gn$ is irreducible so that 
$\Vscr_{tstst}$ is the unique orbital variety
associated to $\Oscr_{tstst}.$ Hence if  $\Vscr_s$
is a union of orbital varieties then
$\ov\Vscr_{s}\supsetneq\Vscr_{tstst}$ 
which is again impossible by proposition \ref{2.3}
since $\tau(s)=\{\be\}$ and 
$\tau(tstst)=\{\al\}.$ 

\subsection{}\label{3.4}
For searching the situation in $D_4$ we use the 
calculations in \cite{Sp}.
Let $s_{\sr 3}$ be the reflection 
giving $s_{\sr 3}(\al_i)=\al_i+\al_{\sr 3}$
for $i=1,2,4$ and $s_{\sr 1},s_{\sr 2},s_{\sr 4}$ the rest 
fundamental reflections. Let us parameterize 
nilpotent orbits in $D_4$
by the partitions corresponding to their 
Jordan form. Then there are 
$\Oscr_1\leftrightarrow (3,3,1,1)$
and $\Oscr_2\leftrightarrow (3,2,2,1)$ such that
$\Oscr_2\st\ov\Oscr_1.$ 
There are only 2 orbital varieties
associated to $\Oscr_2.$
The elementary calculations show that these are
$\Vscr_{s_{\sr 1}s_{\sr 2}s_{\sr 4}s_{\sr 3}s_{\sr 1}s_{\sr 2}s_{\sr 4}}$ and
$\Vscr_{s_{\sr 3}s_{\sr 1}s_{\sr 2}s_{\sr 4}s_{\sr 3}s_{\sr 1}s_{\sr 2}s_{\sr 4}s_{\sr 3}}.$
Note that
$\tau(s_{\sr 1}s_{\sr 2}s_{\sr 4}s_{\sr 3}s_{\sr 1}s_{\sr 2}s_{\sr 4})
=\{\al_{\sr 1},\al_{\sr 2},\al_{\sr 4}\}$ and 
$\tau(s_{\sr 3}s_{\sr 1}s_{\sr 2}s_{\sr 4}s_{\sr 3}s_{\sr 1}s_{\sr 2}s_{\sr 4}s_{\sr 3})=
\{\al_{\sr 3}\}.$
We also have that
$\Vscr_{s_{\sr 1}s_{\sr 3}s_{\sr 1}}$
is associated to $\Oscr_1.$ Note that 
$\tau(s_{\sr 1}s_{\sr 3}s_{\sr 1})=\{\al_1,\al_3\}.$ 

Again, by \ref{3.1} if $\ov\Vscr_{s_{\sr 1}s_{\sr 3}s_{\sr 1}}$ is a union of orbital
varieties it must include at least one of
$\Vscr_{s_{\sr 1}s_{\sr 2}s_{\sr 4}s_{\sr 3}s_{\sr 1}s_{\sr 2}s_{\sr 4}},\ 
\Vscr_{s_{\sr 3}s_{\sr 1}s_{\sr 2}s_{\sr 4}s_{\sr 3}s_{\sr 1}s_{\sr 2}s_{\sr 4}s_{\sr 3}},$
which is impossible by proposition \ref{2.3}.

Again, as in the case of $B_2,$ there is no intermediate nilpotent orbit $\Oscr\pr$
between $\Oscr_1$ and $\Oscr_2.$
\subsection{}\label{3.5}
Now we are ready to show  
\begin{prop} In a semi-simple Lie algebra $\gog$
having a factor not of type $A_n$
there exists an orbital variety such that 
its closure is not a union of orbital varieties.
\end{prop}
\Pf
Our proof is based on the previous computations and proposition
\ref{2.2}.

Indeed,  since orbital variety as well as its closure
in a semisimple Lie algebra is just a direct product
of corresponding simple factors, it is enough 
to prove the proposition for a
simple Lie algebra not of type $A_n$. So,
let $\gog$ be a simple Lie algebra
not of type  $A_n.$
\par
For $\gog$ of type $G_2$ we have shown the existence
of such orbital variety in \ref{3.3}.
\par
If $\gog$ is not of type $G_2,$
then there exist $\Iscr\st \Pi$ of type	
$B_2$ or of type $D_4.$ Let us denote 
simple reflections in $\Iscr$ as in case of $B_2$ in \ref{3.3} if 
$\gog$ is of type
$B_n,\ C_n$ or $F_4,$ and
as in case of $D_4$ in \ref{3.4} otherwise. Let us denote 
$\Oscr_2:=\Oscr_{tst}$ in the case of $B_2$ and keep the notation
$\Oscr_2$ in the case of $D_4$ as in \ref{3.4}.  
Set
$$w_{\sr \Iscr}:=
\begin{cases}
s & {\rm if}\ \gog\ {\rm is\ of\ type}\ B_n,\ C_n,\ {\rm or}\ F_4,\\
                s_{\sr 1}s_{\sr 3}s_{\sr 1} & {\rm otherwise}.\\
\end{cases}
$$
Recall the notion of $D_{\Iscr}$ from \ref{2.2} and set 
$d_m$ to be the maximal element of $D_{\Iscr}.$ 
Such element is unique by the uniqueness of the 
longest element in $W$. 
\par
We will show  that $\ov\Vscr_ {w_{\sr \Iscr}d_m}$ is 
not a union of orbital varieties.

By the construction 
$$\ov\Vscr_{w_{\sr \Iscr}d_m}=
\ov{\bB(\gn_{\sr \Iscr}\cap^{w_{\sr \Iscr}}\gn_{\sr \Iscr})}.\qquad (*)$$
Hence for every $\Vscr$ such that $\pi_{\sr \Iscr}(\ov\Vscr)=
\ov{\Vscr_{w_{\sr \Iscr}}^\Iscr}$ one has
$\ov\Vscr\supset\Vscr_{ w_{\sr \Iscr}d_m}.$

Assume that 
$\ov\Vscr_{w_{\sr \Iscr}d_m}\setminus 
\Vscr_{w_{\sr \Iscr}d_m}=
\cup\Vscr_i$ for some orbital varieties $\Vscr_i.$ By the previous note 
$$\pi_{\sr \Iscr}(\ov \Vscr_i)\subsetneq \ov\Vscr_{w_{\sr \Iscr}}^{\Iscr}.\qquad (**) $$

Now take some point $X\in 
\ov\Vscr_{w_{\sr \Iscr}}^{\Iscr}\cap \Oscr_2.$  Consider it as a point
of $\gn_{\sr \Iscr}\st \gn.$ We denote it by $\hat X$ when we consider it as a point
of $\gog.$ Then 
$\hat X\in\ov\Vscr_{w_{\sr \Iscr}d_m}$ by $(*).$ Moreover,
$\Oscr_{\hat X}\ne\Oscr_{w_{\sr \Iscr}d_m},$ hence, 
$\hat X\in\ov\Vscr_{w_{\sr \Iscr}d_m}\setminus 
\Vscr_{w_{\sr \Iscr}d_m}.$
By our assumption there exist $\Vscr_i$ such that
$\hat X\in \Vscr_i.$ By theorem \ref{2.2} $\pi_{\sr \Iscr}(\ov\Vscr_i)$
is a closure of an orbital variety, and by $(**)$   
$\pi_{\sr \Iscr}(\ov\Vscr_i)\subsetneq \ov{\Vscr_{w_{\sr \Iscr}}^\Iscr}.$
Since $X\in\pi_{\sr \Iscr}(\ov\Vscr_i)$ and there is no intermediate
nilpotent orbits between $\Oscr_{\Vscr_{w_\Iscr}^\Iscr}$ and $\Oscr_2$
(by the notes in the end of case $B_2$ in \ref{3.3} and of \ref{3.4}) we get that
$\pi_{\sr \Iscr}(\ov\Vscr_i)$ is the closure of some orbital variety 
associated to $\Oscr_2.$
This  contradicts our computations in \ref{3.3}, \ref{3.4}.
\QED

\subsection{}\label{3.6}
Now let us consider the situation $\gog=\gs\gl_n.$
Here, as it is shown in \cite[4.1.8]{M1}, one has
\begin{prop}
Let $\gog=\gs\gl_n.$ Then for every two nilpotent orbits $\Oscr_1,\Oscr_2$ such that
$\ov\Oscr_2\subsetneq\ov\Oscr_1$ and for every $\Vscr_1$
associated to $\Oscr_1$ there exist $\Vscr_2$ associated to
$\Oscr_2$ such that $\Vscr_2\st\ov\Vscr_1.$
\end{prop}
Therefore, the argument we use in other cases cannot work for $\gs\gl_n.$
Moreover, modulo this proposition conjecture \ref{1.2} is equivalent
to the equidimensionality of $\ov\Vscr\cap\Oscr$ for any $\Oscr$
in the closure of the nilpotent orbit, $\Vscr$ is associated to. 

As it is shown in \cite[2.3]{M2},
if $\Vscr$ is a Richardson component (that is
its closure is a nilradical of a standard parabolic subgroup)
then its closure is a union of orbital varieties. Note that in our counterexamples
 \ref{3.3}, \ref{3.4} all the orbital varieties in question are Richardson.
This demonstrates again, that the situation in $\gs\gl_n$ is different from other cases.
As well it is shown in \cite[3.15]{M3} that if $\Vscr$
is associated to a nilpotent orbit of nilpotent order 2 then its closure
is a union of orbital varieties. 
As we mentioned in \ref{1.2} these results
together with computations for low ranks support conjecture \ref{1.2}.

\end{document}